\documentclass[11pt]{pjm}
\usepackage{amsmath,amssymb}
\usepackage[latin1]{inputenc}

\newcommand{\ZZ}{\mathbb{Z}}
\newcommand{\RR}{\mathbb{R}}
\newcommand{\CC}{\mathbb{C}}

\newcommand{\Aff}{\mathsf{A}}

\newcommand{\ann}{\operatorname{Ann}}

\newcommand{\witt}{\textsf{W}}

\newtheorem{theo}{Theorem}[section]
\newtheorem{prop}[theo]{Proposition}
\newtheorem{lemma}[theo]{Lemma}

\theoremstyle{definition}

\newtheorem{defi}[theo]{Definition}

\title[Classification of some simple pre-Lie algebras]{Classification
  of some simple graded pre-Lie algebras of growth one}

\date{Received date / Revised version date}
\author{Frédéric Chapoton}
\address{LACIM\\
Université du Québec à Montréal\\
CP 8888 succ. centre ville\\
MONTRÉAL QUÉBEC H3C 3P8\\
CANADA}
\email{chapoton@math.uqam.ca}

\begin{document}
\maketitle
\begin{abstract}
  We classify a class of infinite-dimensional simple graded pre-Lie
  algebras on the graded vector space underlying the algebra of
  Laurent polynomials, with a specific form for the product.
\end{abstract}

\setcounter{section}{-1}

\section{Introduction}

The notion of pre-Lie algebra can be seen as a weakened form of
associative algebra, still giving a Lie algebra by anti-symmetrization.
It has been introduced by Gerstenhaber in his work on deformations of
algebras \cite{gerst}. More recently, it has been studied from the
point of view of operad theory and seen to be related to rooted trees
\cite{rooted}.

The motivating problem for this paper is the classification of
infinite-dim\-ensional simple graded pre-Lie algebras of finite growth
over $\CC$. This seems to be a difficult problem, as it implies in
particular the same classification for associative algebras. The
analogous problem for Lie algebras has been solved by Mathieu in
\cite{mathieu2,mathieu3}. The first step in the case of Lie algebras
was the classification for growth not greater than $1$ done in
\cite{mathieu0,mathieu1}. As a modest first step to the similar
classification problem for pre-Lie algebras, we obtain here the
classification of simple graded pre-Lie algebras such that
\begin{enumerate}
\item The underlying graded vector space is $E=\oplus_{i\in \ZZ} \CC e_i$.
\item The product is given by $ e_i \circ e_j = f(i) g(j) e_{i+j},
  \hfill (\#)$\\ for some complex-valued functions $f$ and $g$ on
  $\ZZ$.
\end{enumerate}
This Ansatz for the product has no special meaning, except that it
allows for a full classification, which is unknown in the general
case.
\begin{theo}
  Any simple graded pre-Lie algebra structure on $E$ of the form
  $(\#)$ is isomorphic either to the algebra $A_a$ defined by
  \begin{equation}
     e_i \circ e_j = (1+a i) e_{i+j}
  \end{equation}
  for some complex number $a$, or to the algebra $B_b$ defined by
  \begin{equation}
    e_i \circ e_j =\frac{i}{1+b j} e_{i+j}
  \end{equation}
  for some complex number $b$ which is not the inverse of an integer.

  All these algebras are pairwise non-isomorphic, with the exceptions
  $A_a \simeq A_{-a}$ and $B_{b}\simeq B_{-b}$.
\end{theo}

One can see that the associated Lie algebras are isomorphic to the
Witt algebra $\witt$, except for $A_0$ with trivial Lie algebra. The
pre-Lie algebras $A_a$ and $B_b$ have been considered before by
Kuperschmidt in his study of a pre-Lie (``quasiassociative'')
structure on the Virasoro algebra \cite{kuper}. There remains the
possibility of other isomorphism classes of simple pre-Lie algebra
structures on $E$ with no presentation fitting our Ansatz.

The first section collects basic general results on pre-Lie algebras.
The second section contains the proof that any simple graded pre-Lie
algebra fitting the Ansatz is isomorphic either to an algebra $A_a$ or
an algebra $B_b$. In the third section, possible isomorphisms between
the algebras $A_a$ and $B_b$ are determined. The last section deals
with injective morphisms into pre-Lie algebras of vector fields on the
line.

\section{General properties of simple pre-Lie algebras}

This section contains basic material on pre-Lie algebras, some of
which can already be found in \cite{burde_reduc,burde_solv}, which are
concerned with finite dimensional pre-Lie algebras (under the name
left-symmetric algebras).

\begin{defi}
  A pre-Lie algebra is a vector space $L$ endowed with $ \circ: L
  \otimes L \to L $ satisfying
  \begin{equation}
    (x \circ y)\circ z-x \circ (y \circ z)=
    (x \circ z)\circ y-x \circ (z \circ y) \quad \forall x,y,z \in
    L.
  \end{equation}
\end{defi}

An \textit{ideal} in a pre-Lie algebra is a two-sided ideal for
$\circ$.

\begin{defi}
  A graded pre-Lie algebra with non-zero product and without any
  non-trivial homogeneous ideal is called simple.
\end{defi}

\begin{lemma}
  In a pre-Lie algebra $L$, the subspaces
  \begin{align}
    \ann(L) &= \{ x \in L \mid \forall y \in L \quad y \circ x = 0\}\\
    \text{ and }\quad L \circ L &=\text{span of } \{ x \circ y \mid
    x,y \in L\}
  \end{align}
  are ideals and are preserved by any automorphism of $L$.
\end{lemma}
\begin{proof}
  These subspaces are clearly invariant by any automorphism. The only
  non-trivial check is that $x\circ y$ is in $\ann(L)$ for
  $x\in\ann(L)$ and $y\in L$. Indeed, one has, for $z \in L$,
  \begin{equation}
  z\circ (x\circ y)=z\circ (y\circ x)-(z\circ y)\circ x+(z\circ x)\circ y,
  \end{equation}
  which vanishes by definition of $\ann(L)$.
\end{proof}

In a graded pre-Lie algebra $L$, $\ann(L)$ and $L \circ L$ are therefore
homogeneous ideals. 

\begin{lemma}
  In a simple graded pre-Lie algebra $L$, $\ann(L)$ is zero and
  $L \circ L$ is $L$.
\end{lemma}
\begin{proof}
  By simplicity, $\ann(L)$ and $L \circ L$ can only be zero or $L$. If
  either $\ann(L)$ is $L$ or $L \circ L$ is zero, then the product is
  zero, which is excluded by definition.
\end{proof}

\begin{lemma}
  \label{desc}
  Let $L$ be a graded pre-Lie algebra. If the associated graded Lie
  algebra is graded simple then $L$ is graded simple as a pre-Lie
  algebra.
\end{lemma}
\begin{proof}
  A non-trivial homogeneous pre-Lie ideal would give a non-trivial
  homogeneous Lie ideal.
\end{proof}

\begin{lemma}
  \label{deri}
  Let $e$ be an \textrm{invariant vector} in a pre-Lie algebra $L$,
  \textit{i.e.} $e\circ x=0$ for all $x\in L$. Then the map $R_e : y
  \mapsto y \circ e$ is a derivation of the pre-Lie product of $L$.
\end{lemma}
\begin{proof}
  Write the pre-Lie axiom for $(x \circ y) \circ e$.
\end{proof}

\section{Reduction to cases $A$ and $B$}

The condition that formula $(\#)$ defines a pre-Lie
algebra can be stated as the vanishing of the numbers
$(C_{i,j,k})_{i,j,k \in \ZZ}$, where $C_{i,j,k}$ is the coefficient of
$e_{i+j+k}$ in
\begin{equation}
  (e_i \circ e_j) \circ e_k - e_i \circ (e_j \circ e_k) 
- (e_i \circ e_k) \circ e_j + e_i \circ (e_k \circ e_j).
\end{equation}
Remark that the numbers $C_{i,j,j}$ are always zero because
$C_{i,j,k}=-C_{i,k,j}$. Explicitly, $C_{i,j,k}$ is given by
\begin{equation}
\label{expanded}
  f(i)\bigg( \big(f(i+j) - f(i+k)\big)g(j)g(k) 
     + \big(f(k)g(j)- f(j)g(k)\big)g(j+k)\bigg).
\end{equation}

One can now use the general properties proved in the previous section
to deduce properties of $f$ and $g$ from the simplicity of $E$.

As $\ann(E)$ is zero, the function $g$ can not vanish anywhere. One
can (and will) in particular assume without restriction that $g(0)=1$
by dividing the basis and the function $g$ by $g(0)$.

As $E \circ E$ is $E$, the function $f$ cannot be identically zero. Let
us discuss according to the value of $f(0)$.

\subsection{Case A: $f(0)\not=0$}

When $f(0)\not=0$, one can assume without restriction that $f(0)=1$ by
dividing the basis and the function $f$ by $f(0)$.

\begin{lemma}
  In this case, one has $g(k)=1$ for all $k$.
\end{lemma}
\begin{proof}
  Write the condition $C_{0,0,k}=0$ and use that $g(k)\not=0$.
\end{proof}

Therefore the vanishing of formula (\ref{expanded}) reduces to
\begin{equation}
  \label{simplified}
  f(i)\bigg(f(i+j)-f(j)-\big(f(i+k)-f(k)\big)\bigg)=0.
\end{equation}

Let us introduce a family of simple graded pre-Lie algebras $A_a$
indexed by a complex number $a$. They are defined on the vector space
$E$ by the formula $ e_i \circ e_j = (1+ a i) e_{i+j}$. One can check
that, for all $a\not=0$, the associated Lie algebra of $A_a$ is the
Witt algebra $\witt$, which is a simple graded Lie algebra. It follows
by Lemma \ref{desc} that $A_a$ is a simple graded pre-Lie algebra for
$a\not=0$. The graded-simplicity of $A_{0}$ is easy to prove and the
associated Lie algebra has zero product.

\begin{prop}
  Any simple graded pre-Lie algebra structure on $E$ with product 
  \begin{equation}
      e_i \circ e_j = f(i) e_{i+j},
  \end{equation}
  for some complex-valued function $f$ on $\ZZ$ with $f(0)=1$, is
  isomorphic to one of the $A_a$.
\end{prop}
\begin{proof}
  Let us first assume that $f(i)=0$ for all $i\not=0$. Then the
  subspace $\oplus_{i\not=0} \CC e_i$ is a non-trivial homogeneous
  ideal of $E_f$, which is a contradiction.
  
  One can therefore assume that there exists $i\not=0$ with
  $f(i)\not=0$. We do the proof when there exists such a positive $i$,
  the negative case being strictly analog. Let $d$ be the smallest
  positive $i$ such that $f(i)\not=0$ and let $f(d)=1+a$ for some
  complex number $a\not=-1$.
  
  If $d=1$, then formula (\ref{simplified}) for $i=1$, $j=0$ and $k$ implies
  that $f(k)=1+k a$ for all $k$. This is the pre-Lie algebra $A_a$
  introduced above.
  
  If $d>1$, then $f(1)=0$ and formula (\ref{simplified}) for $i=d$, $j=0$ and
  $k=1$ implies that $f(d+1)=a$. If $a=0$ then $f$ is periodic of
  period $d$ and therefore
  \begin{equation}
    \begin{cases}
      f(i)=1 \text{  if  } d \mid i,\\
      f(i)=0 \text{  otherwise.}
    \end{cases}
  \end{equation}
  In this case, the subspace $\oplus_{i\not \equiv 0 [d]} \CC e_i$ is a
  non-trivial homogeneous ideal of $E_f$, which is a contradiction.
  
  One can therefore assume that $d>1$ and $a\not=0$. If $d=2$, the
  sequence of values of $f$ on non-negative integers starts with
  \begin{equation}
    1,0,1+a,a,1+2a,\dots,
  \end{equation}
  and $a\not=0$ implies $1+2a=a-1$, \textit{i.e.} $a=-2$ and
  $f(k)=1-k$ for all $k$. This is the pre-Lie algebra $A_{-1}$.
  
  The only unsettled case is when $d\geq 3$ and $a\not=0$.
  Then the sequence of values of $f$ on non-negative integers starts
  with
  \begin{equation}
    1,\underbrace{0,0,\dots}_{d-1 \text{ zeroes}},1+a,a,a,\dots,
  \end{equation}
  and $a\not=0$ implies $a=a-1$, which is a contradiction.
\end{proof}

\subsection{Case B: $f(0)=0$}

When $f(0)=0$, there exists $i\not=0$ such that $f(i)\not=0$. The
vector $e_0$ is invariant, so it gives a derivation by Lemma
\ref{deri}. The vanishing of $C_{i,j,0}$ is then equivalent to the
following equation:
\begin{equation}
  \label{quasiperiod}
  f(i)\bigg(f(i+j)-f(i)-f(j)\bigg)=0.
\end{equation}

Let us introduce a family of simple graded pre-Lie algebras $B_b$ for
a complex number $b$ which is not the inverse of an integer. They are
defined on the vector space $E$ by the formula $ e_i \circ e_j = i/(1+
b j) e_{i+j}$. By the change of basis $\bar{e}_i=(1+ b i) e_i$, this
is equivalent to
\begin{equation}
  \label{modif}
  \bar{e}_i \circ \bar{e}_j =\frac{ i(1+b i)}{1+
  b(i+j)} \bar{e}_{i+j},
\end{equation}
which is the formula given by Kuperschmidt \cite{kuper}.

One can check, by using formula (\ref{modif}), that the associated
Lie algebra of $B_b$ is the Witt algebra $\witt$ for all $b$. It
follows by Lemma \ref{desc} that $B_b$ is a simple graded pre-Lie
algebra for any admissible $b$.

\begin{prop}
  Any simple graded pre-Lie algebra structure on $E$ with product
  \begin{equation}
      e_i \circ e_j = f(i)g(j) e_{i+j},
  \end{equation}
  for some complex-valued functions $f,g$ on $\ZZ$ with $f(0)=0$, is
  isomorphic to one of the $B_b$.
\end{prop}
\begin{proof}
  
  There exists $i\not=0$ such that $f(i)\not=0$. Assume for example
  that there exists a positive such $i$ and let $d$ be the smallest
  positive integer such that $f(d)\not=0$. Then one can assume without
  restriction that $f(d)=d$, by dividing the basis and $f$ by
  $f(d)/d$. It follows from (\ref{quasiperiod}) that $f(i)-i$ is a
  periodic function of $i$ with period $d$. It has also period $d+1$,
  because $f(d+1)$ is either $d$ or $2d$, which is non-zero. This
  implies that $f(i)-i$ is constant. As $f(0)=0$, one has $f(i)=i$ for
  all $i$.

  
  Let us define a complex-valued function on $\ZZ$ by
  $h(i)=1/g(i)$. The vanishing of (\ref{expanded}) becomes 
  \begin{equation}
    k\big(h(j+k)-h(k)\big)=j\big(h(j+k)-h(j)\big).
  \end{equation}
  As a special case, $h(j)+h(-j)=2$ for all $j$. Let $h(1)=1+b$ for
  some complex number $b$. Then $h(-1)=1-b$. One proves by recursion
  that $h(j)=1+b j$ for all $j$.

  To summarize, we have obtained the following law:
  \begin{equation}
    e_i \circ e_j = \frac{i}{1+b j} e_{i+j},
  \end{equation}
  which is defined only if $b$ is not the inverse of an integer, and
  is the pre-Lie algebra $B_b$ introduced above.

\end{proof}

\section{Isomorphism types}

In this section, the algebras $A_a$ and $B_b$ are not considered as
graded, \textit{i.e.} one works in the category of not necessarily
graded pre-Lie algebras.

\smallskip

The map $e_i \mapsto e_{-i}$ defines an isomorphism between $A_{a}$
and $A_{-a}$. The map $e_i \mapsto -e_{-i}$ defines an isomorphism
between $B_{b}$ and $B_{-b}$.

A useful remark is that $A_0$ is the only associative algebra among
the $A_a$ and $B_b$. In fact, $A_0$ is isomorphic to the commutative
algebra $\CC[x,x^{-1}]$ of Laurent polynomials.

\begin{prop}
  \label{isomo}
  The only non-trivial (not necessarily graded) isomorphisms between
  the $A_a$ and $B_b$ pre-Lie algebras are the above isomorphisms $A_a
  \simeq A_{-a}$ and $B_{b}\simeq B_{-b}$.
\end{prop}

Before proving the proposition, we state a crucial lemma. Let $L$ be a
pre-Lie algebra which is one of the $A_a$ or $B_b$.

\begin{lemma}
  \label{crux}
  Let $v\in L$. If $R_v: x\mapsto x \circ v$ is a locally finite
  endomorphism of $L$, then $v$ is proportional to $e_0$.
\end{lemma}
\begin{proof}
  One can assume that $v\not=0$. Let $v=\sum_{i=k}^{\ell} \lambda_i
  e_i$ with $\lambda_{k}$ and $\lambda_{\ell}$ non-zero. Let us show
  that $\ell>0$ implies that $R_v$ is not locally finite.
  
  Let $N$ be a large positive integer. Consider the sequence
  \begin{equation}
    e_{N},R_v(e_N),R_v^2(e_N),\dots
  \end{equation}
  By recursion, one can see that for all $p\geq 0$, $R_v^p(e_N)$ is a
  finite sum of the form $\sum_{j\leq N+p \ell}\mu_j e_{j}$ with
  $\mu_{N+p \ell}\not=0$, provided that $N$ is larger than $-1/a$ in
  the case of the $A_a$ algebra. Therefore these vectors are linearly
  independent. This is a contradiction with local finiteness of $R_v$.
  Hence we can assume that $\ell \leq 0$. By the same proof with
  reverse grading, one can assume that $k\geq 0$, so that $v\in\CC
  e_0$.
\end{proof}

Let us go back to the proof of proposition \ref{isomo}.

\begin{proof}
  Let $L$ be a pre-Lie algebra which is one of the $A_a$ or $B_b$. By
  Lemma \ref{crux}, the subspace $\CC e_0$ is invariant by any
  automorphism. Let $v$ be a non-zero multiple of $e_0$. We now
  discuss according to the vanishing of $v\circ v$.
  
  If $v \circ v =0$, then the algebra must be one of the $B_b$. The
  eigenspace decomposition of $R_v$ gives back the homogeneous
  decomposition, up to reversal of the grading. Choose one of the two
  opposite gradings and let $w\not=0$ be in $\CC e_{1}$. Then $(w
  \circ w)\circ w$ and $w \circ (w \circ w)$ are two non-zero vectors
  in $\CC e_{3}$. Their proportionality coefficient gives back $b$ or
  $-b$, according to the choice of grading.
  
  If $v \circ v \not=0$, then the algebra must be one of the $A_a$.
  One can further assume that $a\not=0$, because $A_0$ is the only
  associative algebra among the $A_a$. Then $e_0$ is the unique vector
  in $\CC e_0$ such that $e_0 \circ e_0=e_0$. The eigenvalues of $R_{e_0}$
  give back the sequence $1+ a i$ up to a choice of grading. So we
  recover $a$ up to sign.
\end{proof}

\section{Vector fields realizations}



We determine here all realizations of the algebras $A_a$ and $B_b$ as
subalgebras of the pre-Lie algebra $V_1(\Omega)$ of analytic complex
vector fields on an connected open subset $\Omega$ of the complex
affine line $\Aff^1(\CC)$.

Let us start with the $A$ case:

\begin{prop}
  The algebra $A_a$ for $a\not=0$ is isomorphic to the subalgebra of
  $V_1(\CC\setminus \RR^-)$ with basis
  \begin{equation}
    e_i= x^{1+a i} \frac{\partial}{\partial x}.
  \end{equation}
  Up to an affine change of variable, this is the unique realization of
  $A_a$ as a pre-Lie algebra of vector fields on the line.
  
  There is no injection of $A_0$ into the pre-Lie algebra of analytic
  vector fields on an open subset of $\Aff^1(\CC)$.
\end{prop}
\begin{proof}
  Let $\psi$ be a map from $A_a$ to such a pre-Lie algebra. Because of
  $e_0 \circ e_0 =e_0$, $\psi(e_0)$ must belong to $\CC \partial_x
  \oplus \CC x\partial_x$ with non-vanishing component in $\CC
  x\partial_x$. One can assume by a translation of the variable $x$
  that $\psi(e_0)$ belongs to $\CC x\partial_x$. Then, in fact,
  $\psi(e_0)=x\partial_x$. Because $e_i \circ e_0 =(1+a i) e_i$,
  $\psi(e_i)$ is in $\CC x^{1+a i}\partial_x$ for all $i$. Hence one
  must restrict to a simply-connected open subset of $\CC^*$, and
  injectivity is only possible if $a\not=0$. By a dilatation of the
  variable $x$, one can assume that $\psi(e_1)=x^{1+a}\partial_x$.
  Then is is easy to see that necessarily $\psi(e_i)=x^{1+a
    i}\partial_x$ for all $i$. Finally, it is easy to check that
  $\psi$ is a pre-Lie algebra map.
\end{proof}

Now, consider the $B$ case:

\begin{prop}
  The algebra $B_0$ is isomorphic to the subalgebra of $V_1(\CC)$ with
  basis
  \begin{equation}
    e_i=\exp(i x) \frac{\partial}{\partial x}.
  \end{equation}
  Up to an affine change of variable, this is the unique realization of
  $B_0$ as a pre-Lie algebra of vector fields on the line.
  
  There is no injection of $B_b$ into the pre-Lie algebra of analytic
  vector fields on an open subset of $\Aff^1(\CC)$ for $b\not=0$.
\end{prop}
\begin{proof}
  Let $\psi$ be a map from $B_b$ to such a pre-Lie algebra. Because of
  $e_0 \circ e_0 =0$, $\psi(e_0)$ must be in $\CC \partial_x$. By a
  dilatation of the variable $x$, one can assume that
  $\psi(e_0)=\partial_x$. Then $e_i \circ e_0=i e_i$ implies that
  $\psi(e_i)$ is in $\CC \exp(i x)\partial_x$ for all $i$. By a
  translation of the variable $x$, one can assume that
  $\psi(e_1)=\exp(x)\partial_x$. Then one can check that $\psi(e_1)
  \circ \psi(e_{-1})=-\psi(e_{-1}) \circ \psi(e_{1})$, which implies
  $b=0$. It is easy then to see that $\psi(e_i)=\exp(i x)\partial_x$
  for all $i$. This indeed defines a pre-Lie algebra map.
\end{proof}


\bibliographystyle{alpha}
\bibliography{graded}

\end{document}